\documentclass[12pt,twoside]{article}
\usepackage{amsmath}
\usepackage{amsfonts}
\usepackage{amsthm}
\usepackage{amssymb}
\usepackage{graphicx}
\usepackage{multicol}

\begin{document}
\title{{\normalsize 
{\bf Orthogonal 2-sphere basis of stable 4-sphere}}} 
\author{{\footnotesize Akio Kawauchi}\\ 
\date{} 
{\footnotesize{\it Osaka Central Advanced Mathematical Institute, Osaka Metropolitan University}}\\  
{\footnotesize{\it Sugimoto, Sumiyoshi-ku, Osaka 558-8585, Japan}}\\  
{\footnotesize{\it kawauchi@omu.ac.jp}}} 
\date\,  
\maketitle 
\vspace{0.1in} 
\baselineskip=9pt 
\newtheorem{Theorem}{Theorem}[section] 
\newtheorem{Conjecture}[Theorem]{Conjecture} 
\newtheorem{Lemma}[Theorem]{Lemma} 
\newtheorem{Sublemma}[Theorem]{Sublemma} 
\newtheorem{Proposition}[Theorem]{Proposition} 
\newtheorem{Corollary}[Theorem]{Corollary} 
\newtheorem{Claim}[Theorem]{Claim} 
\newtheorem{Definition}[Theorem]{Definition} 
\newtheorem{Example}[Theorem]{Example} 

\begin{abstract} Every stable 4-sphere is identified with the double branched covering space of a trivial surface-knot space. As a result of Wall, it is known that any two orthogonal bases of every stable 4-sphere are transformed into each other by an orientation-preserving diffeomorphism of the stable 4-sphere. In this paper another proof of Wall's result is presented, strengthened in the sense that the lift of an equivalence of the trivial surface-knot space can be taken as the diffeomorphism. Two applications are made. The first shows that every orientation-preserving diffeomorphism of every stable 4-sphere is nothing but the double branched covering lift of an equivalence of a trivial surface-knot space up to a smooth isotopy and a composition with an identity-shift. The second gives a similar result for TOP stable 4-spheres. Here, even if it is a smooth 4-manifold, unless it is diffeomorphic to the stable 4-sphere, the TOP trivial surface-knot space cannot be smooth.

\phantom{x}

\noindent{\footnotesize{\it Keywords:} Stable 4-sphere,\, 
Trivial surface-knot,\, O2-handle basis,\, Branched covering,\, Orthogonal basis,\,  O2-sphere basis.} 

\noindent{\it Mathematics Subject Classification 2020}:  57K45, 57K40
\end{abstract}

\baselineskip=15pt

\bigskip

\noindent{\bf 1. Introduction}

\begin{figure}[hbtp]
	\begin{center}\includegraphics[width=7.5 cm,height=5.5 cm]{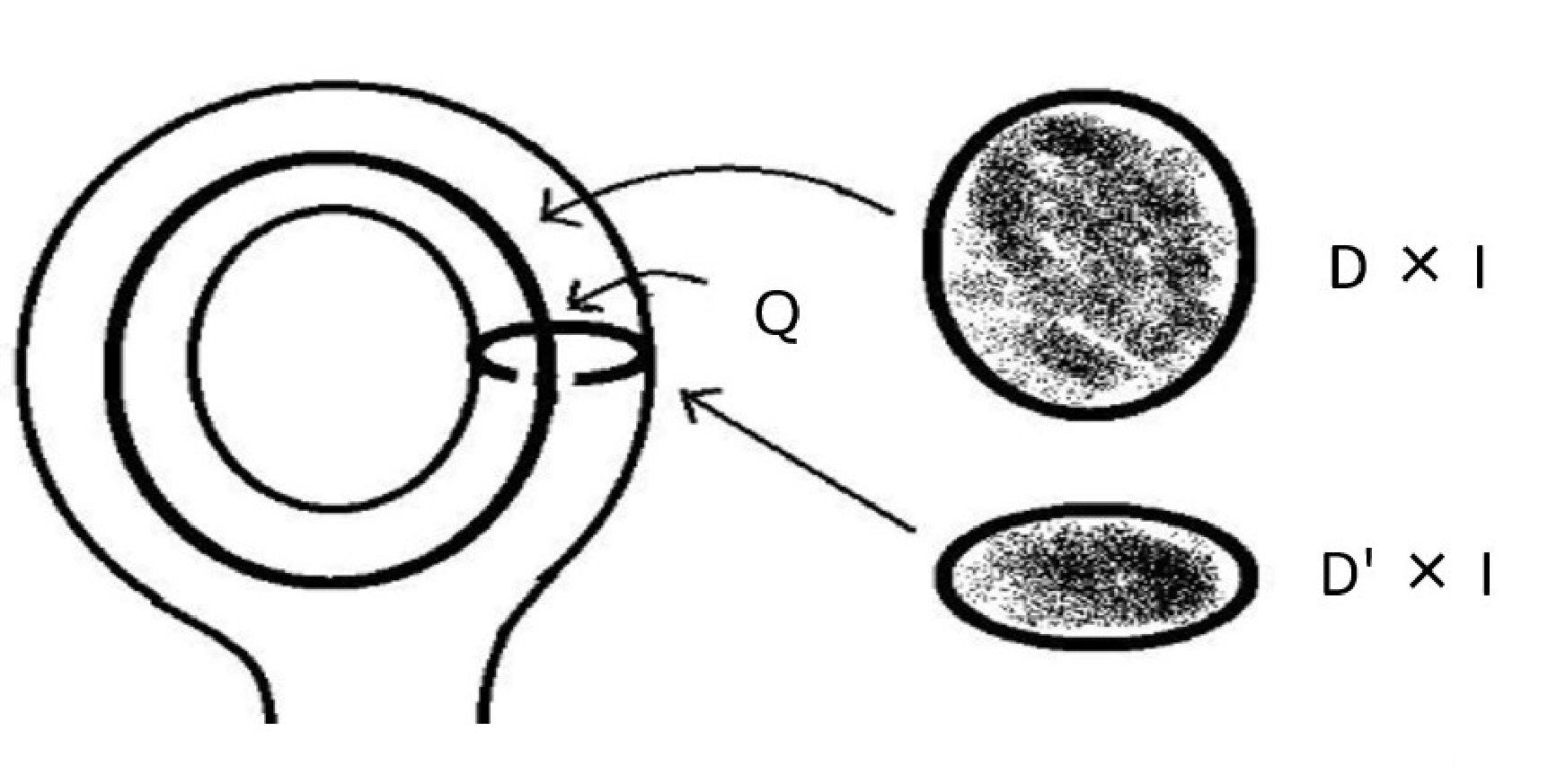}\end{center}
	\caption{O2-handle pair $(D\times I, D'\times I)$}
	\label{fig1}
\end{figure}

The {\it stable 4-sphere of genus} $n$ is  
the  connected sum 4-manifold $\Sigma=\Sigma(n)$ of $n$ copies of 
$S^2\times S^2$. The stable 4-space $\Sigma(n)$
 is canonically diffeomorphic to the double branched covering space $S^4(F)_2$ of the 4-sphere $S^4$ branched along a trivial surface-knot $F$ of genus $n$ in $S^4$, \cite{1}.  The pair $(S^4,F)$ is called a {\it trivial surface-knot space} 
of genus $n$.
Here, a {\it trivial surface-knot of genus} $n$ in $S^4$ is a surface-knot 
bounding a handlebody of genus $n$ smoothly embedded in $S^4$, whose pair 
$(S^4, F)$ is called a {trivial surface-knot space of genus} $n$. 
An {\it orthogonal basis} of the stable 4-sphere $\Sigma(n)$ 
 is a pairwise basis $(x_*,x'_*) = \{(x_i,x'_i)|i = 1, 2, \dots, n\}$ of the second 
integral homology group $H_2(\Sigma(n); Z)$ which is a free abelian group of rank 
$2n$ such that the intersection numbers in $\Sigma(n)$ have 
$\mbox{Int}(x_i, x_j) = \mbox{Int}(x_i, x'_j) =  
\mbox{Int}(x'_i, x_j) = \mbox{Int}(x'_i, x'_j) = 0$ for all $i, j$ except for that 
$\mbox{Int}(x_i, x'_i) = \mbox{Int}(x'_i, x_i) = 1$ for all $i$.  
It is known as Wall's result that any two orthogonal bases of every stable 
4-sphere $\Sigma$ are transformed into each other by an orientation-preserving 
diffeomorphism $\sigma$ of the stable 4-sphere $\Sigma$, \cite{2}. 
The {\it standard O2-sphere basis} of 
the stable 4-sphere $\Sigma(n)$ is the 2-sphere pair system 
$(S^2\times1_*, 1 \times S^2_*) = \{(S^2\times1_i, 1 \times S^2_i) 
\vert i = 1, 2, \dots, n\}$ 
in $\Sigma(n)$. An {\it O2-sphere basis} of 
$\Sigma(n)$ is a 2-sphere pair system 
$(S_*, S'_*)= \{(S_i, S'_i) \vert i = 1, 2, \dots, n\}$ of $\Sigma(n)$, sent to the standard O2-sphere basis 
$(S^2 \times 1_*,1 \times S^2_*)$ by 
an orientation-preserving diffeomorphism of $\Sigma(n)$, or equivalently which can 
rewrite the stable 4-sphere $\Sigma(n)$ as the connected sum of the genus one stable 4-spheres $S_i \times S'_i \, (i = 1, 2, \dots, n)$. Wall's result is essentially equivalent to saying that every orthogonal basis $(x_*, x'_*)$ of $\Sigma(n)$ is represented by an O2-sphere basis $(S_*, S'_*)$ of $\Sigma(n)$, i.e., 
$x_i = [S_i]$ and $x'_i = [S'_i]\, (i = 1, 2, \dots, n)$, because any two O2-sphere bases of $\Sigma$ are transformed into each other by an orientation-preserving 
diffeomorphism of $\Sigma(n)$.
In this paper, another proof of Wall's result is presented, strengthened in 
the sense that the diffeomorphism $\sigma$  is taken to be the lift $f'$ of an equivalence 
 $f$  of $(S^4, F)$ to the stable 4-sphere  $S^4(F)_2 = \Sigma(n)$. 
 Here, an {\it equivalence} of $(S^4, F)$ is  
 an orientation-preserving diffeomorphism of $S^4$ keeping the oriented trivial 
surface-knot $F$  set-wise fixed. To state the main theorem (Theorem 1.1), the 
 notion of an orthogonal 2-handle pair, or briefly an O2-handle pair on $F$  is 
needed. A {\it 2-handle} on a trivial surface-knot $F$  is a 2-handle $D \times  I $  on $F$ smoothly embedded in $S^4$  such that 
$(D \times  I) \cap F = (\partial D) \times  I$ for a closed  interval $I$ with $0$ as the center and $D\times 0$ is called the {\it core disk} of the 2-handle $D\times I$  
and identified with D. An {\it orthogonal 2-handle pair}, or briefly an 
{\it O2-handle pair} 
 on $F$  in $S^4$  is a pair 
 $(D\times I, D'\times I)$ of 2-handles $D\times I$ and $D'\times I$ on $F$  which meet $F$  only 
with the  attaching annuli $(\partial D)\times I$ and $(\partial D')\times I$ so that the loops $\partial D$ and $\partial D'$ meet transversely at just one-point $x_0$ and the intersection $(\partial D) \times  I \cap  (\partial D') \times  I$ is diffeomorphic to the square $Q =  \{x_0\} \times I \times I$, Fig.~\ref{fig1}, \cite{3}. For a trivial surface-knot $F$  
of genus $n$ in  $S^4$, an O2-handle basis of $F$  is a system 
$(D_* \times  I, D'_*\times  I)$ of mutually 
disjoint O2-handle pairs $(D_i \times  I, D'_i\times  I) (i = 1, 2, \dots, n)$ on $F$  in $S^4$.  Let $p: \Sigma(n) =  S^4(F)_2   \to  S^4$  be the double branched covering projection branched along $F$, and 
$\alpha$  the nontrivial covering involution of  $S^4(F)_2$.  The preimage $p^{-1}(F)$ of $F$  is the 
fixed point set of $\alpha$, diffeomorphic to $F$  and written by the same notation 
as $F$. 
For every O2-handle basis $(D_* \times  I, D'_*\times  I)$ on $F$  in  $S^4$, the O2-sphere basis $(S(D_*), S(D'_*))$ in $\Sigma(n)$ is constructed so that 
$S(D_i) = D_i \cup \alpha D_i$ and $S(D'_i) = D'_i \cup \alpha D'_i\, (i = 1, 2, \dots, n)$ are the preimages of the disks $D_i$ and $D'_i\, (i = 1, 2, \dots, n)$ in $S^4$  by the double branched covering projection p, respectively, [1]. The orientations of 
 $S(D_i)$ and $S(D'_i)$ are taken with the orientations of $D_i$ and $D'_i$ and the opposite 
orientations of $\alpha D_i$ and $\alpha D'_i$, respectively. The main theorem (Theorem 1.1) is stated as follows.

\phantom{x}

\noindent{\bf Theorem~1.1.} Every stable 4-sphere $\Sigma(n)$ is the  double branched covering space   $S^4(F)_2$ of $S^4$ branched along a trivial surface-knot $F$ of genus $n$. Every orthogonal basis  $(x_* , x'_*)$ of $\Sigma(n)$ is represented by an O2-sphere basis $(S(E_*^x), S({E'}_*^x))$ of $\Sigma(n)$ constructed from an O2-handle basis 
$(E_*^x\times  I, {E'}_*^x \times  I )$  on $F$  in $S^4$. For any two orthogonal bases  $(x_* , x'_*)$ and $(y_* , y'_* )$ of $\Sigma(n)$, there are O2-handle bases 
$(E_*^x\times  I,  {E'}_*^x \times  I) and (E_*^y\times  I, {E'}_*^y \times  I )$ on $F$  in $S^4$  which are transformed into each other by an 
 equivalence $f$   of  $(S^4, F)$ so that  $(x_* , x'_*)$ and $(y_* , y'_* )$ are represented by the O2-sphere bases $(S(E_*^x), S({E'}_*^x))$ and 
 $(S(E_*^y), S({E'}_*^y))$, respectively which are transformed into  each other by the lift $f'$ of  $f$   to   $S^4(F)_2 =\Sigma(n)$. 
 
\phantom{x}

A key to showing Theorem 1.1 is the following lemma.  

\phantom{x}

\noindent{\bf Lemma 1.2 (Key Lemma).} 
For every orthogonal basis  $(x_* , x'_*)$ of the stable 4-sphere $\Sigma(n)$, there is an O2-handle basis $(E_*\times  I, E'_* \times  I )$ on a trivial surface-knot $F$  of genus  $n$ in $S^4$  with $\Sigma(n)=S^4(F)_2$   such that  $(x_* , x'_*)$ is represented by the O2-sphere basis $(S(E_*), S(E'_*))$  of $\Sigma(n)$.  

\phantom{x}

A loop basis of $F$  is a pair system $(e_*, e'_*)$ of oriented simple loop pairs 
$(e_i, e'_i)\, (i =  1, 2, \dots, n)$  on $F$  which represents a basis for $H_1(F; Z)$  such that $e_i \cap  e_j = e'_i\cap  e'_j = e_i \cap  e'_j = \emptyset$  for all distinct 
$i, j$ and $e_i \cap  e'_i$ is one point with the intersection number 
$\mbox{Int}(e_i, e'_i) = +1$ in $F$  for all $i$. Every oriented loop c on $F$  bounds an immersed  surface $C$ in $S^4$  with $C \cap  F= c$. The map  
$q: H_1(F; Z) \to  Z_2$ sending the homology  class $[c]$ in $H_1(F; Z)$ to the 
$Z_2$-self-intersection number of $C$ in $S^4$  with respect to  the $F$-framing is called the $Z_2$-{\it quadratic function} associated with the surface-knot $F$  in $S^4$.  The identity $q(x+y)=q(x)+q(y)+ \mbox{Int}F(x, y)_2$ for all $x, y$ in 
$H_1(F; Z)$ is used for calculation, where $\mbox{Int}F(x, y)_2$ denotes the 
$Z_2$-intersection number of $x$ and $y$ in $F$. A  simple loop basis $(e_*, e'_*)$ of $F$  is {\it spin} if $q(e_i) = q(e'_i)=0$ for all $i$. The loop pair system  
$(\partial D_*, \partial D'_*)$ for any O2-handle basis 
$(D_* \times  I, D'_*\times  I)$ on $F$  is a spin loop basis of $F$. Every handlebody smoothly embedded in $S^4$  is smoothly isotopic to a standard  handlebody in the equatorial 3-sphere $S^3$ of $S^4$. Thus, a trivial surface-knot $F$  in $S^4$ may be taken one in the standard position in $S^3$, where a standard O2-handle basis and a standard spin loop basis on $F$  is taken. Any two spin loop bases on $F$  are transformed into each other by an equivalence of  $(S^4, F)$, \cite[(2.5.1), (2.5.2)]{3}, \cite{5}. Hence every spin loop basis on $F$  bounds the core disk pair system of an O2-handle basis on $F$. The core disk pair systems of any two O2-handle bases on $F$  bounded by the same loop basis are transformed into each other by an  equivalence of  $(S^4, F)$, \cite{6}. Under these preliminaries, the proof of Theorem 1.1  assuming Lemma 1.2 is done as follows.

\phantom{x}

\noindent{\bf Proof of Theorem~1.1 assuming Lemma~1.2.} 
The first claim is shown, \cite{1}. The  second claim is shown by Lemma 1.2. If the third claim is proven, then the  conclusion will be shown. By Lemma 1.2, let 
$x_i = [S(E_i^x)]$, $x'_i = [S({E'}_i^x)]$ and $y_i =  [S(E_i^y)]$, $y'_i = [S({E'}_i^y)]$ for O2-handle bases  $(E_*^x \times  I, {E'}_*^x \times  I)$ and 
$(E_*^y \times  I, {E'}_*^y \times  I)$   on a trivial surface-knot $F$  of genus $n$ in $S^4$. Then there is an equivalence $g_0$ of $(S^4,  F)$ sending the spin loop basis $(\partial E_*^x, \partial {E'}_*^x)$ of $F$  to the spin loop basis 
$(\partial E_*^y, \partial {E'}_*^y)$ of $F$, \cite{3}, \cite{5}. For the O2-handle bases $(g_0E_*^x \times  I, g_0{E'}_*^x \times  I)$ and $(E_*^y \times  I, {E'}_*^y \times  I)$  on  $F$  with the same attaching part in $F$, there is an equivalence $g$ of $(S^4, F)$ such that  $(gg_0E_*^x \times  I, gg_0{E'}_*^x \times  I) = (E_*^y \times  I, {E'}_*^y \times  I)$, \cite{6}. The composite equivalence $f  = gg_0$ of  $(S^4, F)$ lifts to a diffeomorphism $f'$ of $\Sigma(n)= S^4(F)_2$   sending the O2-sphere basis $(S(E_*^x), S({E'}_*^x))$ to the O2-sphere basis $(S(E_*^y), S({E'}_*^y))$ and hence sending $(x_i  , x'_i) = ([S(Eix)], [S({E'}_i^x)])$ to $(y_i, y'_i) = ([S(E_i^y)], [S({E'}_i^y)])$ for all $i$. This completes the  proof of Theorem 1.1 assuming Lemma 1.2. 
                
\phantom{x}                
 
Theorem 1.1 is a revised version of an incorrect claim \cite[Lemma 3.1]{4} that was 
originally intended for use in the paper \cite{1}, which was written instead by using  Wall's result, \cite{2}. 
               
The proof of Key Lemma (Lemma 2.1) is done in Section~2.  Two applications 
 of Theorem 1.1 are done in Sections~3 and 4. In Section~3, it is shown that every orientation-preserving diffeomorphism of $\Sigma$ is nothing but the lift of an equivalence of a trivial surface-knot space $(S^4, F)$ to the double branched covering space $\Sigma=S^4(F)_2$ up to a smooth isotopy and a composition of an identity-shift. In Section~4, a TOP version of Theorem~1.1 is done for every TOP stable 4-sphere of genus $n$ (i.e., topological 4-manifold homeomorphic to the stable 4-space of 
 genus $n$),  Here, even if it is a smooth 4-manifold, unless it is diffeomorphic to the stable 4-sphere of genus $n$, the TOP trivial surface-knot space cannot be smooth.

\phantom{x}

\noindent{\bf 2. Proof of Key Lemma (Lemma~1.2)}

For a disk $D$, let $D^o=D\setminus \partial D$. 
The following lemma gives a basic information on the intersection numbers of 
the lifting O2-sphere bases of  two O2-handle bases of $F$ in $S^4$,  which 
corrects a computation error of \cite[Lemma~3.1]{4}. 

\phantom{x}

\noindent{\bf Lemma~2.1.} For an O2-handle basis $(D_* \times  I, D'_*\times  I)$ on a trivial surface-knot $F$  of genus $n$ in  $S^4$, let
$(k_*,k'_*)=(\partial D_*, \partial D'_*)$ be  the  spin loop basis of $F$.
For a 2-handle $E\times I $ on $F$ in $S^4$, assume that the homology class 
$[e]$ of the simple loop $e=\partial E$ in $F$ 
 is given by the intersection numbers $\mbox{Int}([e],[k'_j]) = s_j$ and $\mbox{Int}([e],[k_j]) = s'_j$ in $F$ 
for some integers  $s_j, s'_j\, (j=1,2,\dots,n)$. 
Then the homology class $[S(E)]$ in $\Sigma$ is written as 
$[S(E)]= \sum_{j=1}^n  (s_j+2m_j)[S(D_j)] +\sum_{j=1}^n (s'_j+2m'_j)[S(D'_j) ]$,
where $m_j$ and $m'_j$ are integers given by the intersection numbers 
$\mbox{Int}(E^o,{D'}_j^o)$ and $\mbox{Int}(E^o,D_j^o)$, respectively.

\phantom{x}
 
\noindent{\bf Proof of Lemma~2.1.} Let $N(F)^c=\mbox{cl}(S^4\setminus N(F))$ 
for a regular neighborhood $N(F)$ of $F$ in $S^4$.  
Consider that the disk  
$E^e=E\cap N(F)^c$ transversely meets  the disks  
$D^e_j=D_j\cap N(F)^c$ and ${D'}_j^o=D'_j\cap N(F)^c$ 
 with the intersection points $E^o\cap D_j^o$ and  $E^o \cap {D'}_j^o$ 
 for all $i$, respectively. 
The intersection number of the lift of $E^e$ and the lift of 
$D^e_j$ to  $S^4(F)_2=\Sigma$ is equal to 
$\mbox{Int}(E^e,D^e_j)+\mbox{Int}(\alpha E^e,\alpha D^e_j)=2\mbox{Int}(E^e,D^e_j)
=2\mbox{Int}(E^o,D_j^o)=2m'_j$. 
Similarly, the intersection number of the lift of $E^e$ and the lift of ${D'}_j^o$ to  
$S^4(F)_2=\Sigma$ is equal to 
$\mbox{Int}(E^e,{D'}_j^o)+\mbox{Int}(\alpha E^e,\alpha {D'}^e_j))
=2\mbox{Int}(E^e,{D'}_j^o)=2\mbox{Int}(E^o,{D'}_j^o)=2m_j$.
By using the intersection numbers $\mbox{Int}([e],[k_j])=s'_j$, 
 $\mbox{Int}([e],[k'_j])=s_j$ in $F$ and examining 
the geometric intersections between 
the lift of $E\cap N(F)$ and the lifts of $D_j\cap N(F)$, $D'_j\cap N(F)$ to  
$S^4(F)_2=\Sigma$, the identities  
$\mbox{Int}([S(E)],[S(D_j)])=s'_j+2m'_j$,  
$\mbox{Int}([S(E)],[S(D'_j)])=s_j+2m_j$ 
are obtained. 
Since $([S(D_*)],[S(D'_*)])$ is an orthogonal basis of $\Sigma$, 
the desired identity is obtained. 
This completes the proof of Lemma~2.1.

\phantom{x}

By using Lemma~2.1, the following lemma is obtained. 

\phantom{x}

\noindent{\bf Lemma~2.2.} For every orthogonal basis  $(x_*,x'_*)$  of $\Sigma$, 
there is an O2-handle basis  $(E_*\times I , E'_*\times I )$ on the trivial  surface-knot $F$ of genus $n$ in $S^4$ such that 
$x_i= [S(E_i)]+ 2A_i+2A'_i$,  $x'_i=[S(E'_i)]+2B_i+2B'_i$ for all $i$, 
where 
$A_i=\sum_{j=1}^n a_{ij}[S(E_j)]$,  $A'_i=\sum_{j=1}^n a'_{ij}[S(E'_j)]$, 
 $B_i=\sum_{j=1}^n b_{ij}[S(E_j)]$, $B'_i=\sum_{j=1}^n b'_{ij}[S(E'_j)]$ 
with some integers $a_{ij}$, $a'_{ij}$, $b_{ij}$, $b'_{ij}$  for all $i$. 

\phantom{x}

\noindent{\bf Proof of Lemma~2.2.} For an O2-handle basis 
$(D_* \times  I, D'_*\times  I)$ on a trivial surface-knot $F$  of genus $n$ in  $S^4$, write $x_i$ and $x'_i$ as integral linear combinations of the homology classes 
$[S(D_*)]$ and $[S(D'_*)]$ such that 
$x_i=\sum_{j=1}^n c_{ij}[S(D_j)]+\sum_{j=1}^n c'_{ij}[S(D'_j)]$,  
$x'_i=\sum_{j=1}^n d_{ij}[S(D_j)]+\sum_{j=1}^n d'_{ij}[S(D'_j)]$
for all $i$.
Since  $(x_* , x'_*)$ and $([S(D_*)], [S(D'_*)])$ are orthogonal bases of $\Sigma$, 
the identities
$\sum_{j=1}^n c_{ij}c'_{ij}=\sum_{j=1}^n d_{ij}d'_{ij}=0$, 
$\sum_{j=1}^n (c_{ij}d'_{ij}+ c'_{ij}d_{ij})=1$ hold 
for all $i$.
Let $(k_*,k'_*)=(\partial D_*, \partial D'_*)$ be  a  spin loop basis of $F$.
Then there is a spin loop basis $(\ell_*,\ell'_*)$ on $F$ such that 
$[\ell_i]=\sum_{j=1}^n c_{ij}[k_j]+\sum_{j=1}^n c'_{ij}[k'_j]$,  
$[\ell'_i]=-\sum_{j=1}^n d_{ij}[k_j]+\sum_{j=1}^n d'_{ij}[k'_j]$ for all $i$  
in $H_1(F;Z)$. In fact,  the homology classes $z_i$, $z'_i$ given by the right hands of 
these identities, respectively form a symplectic basis $(z_*,z'_*)$ of $H_1(F;Z)$. Then  
a loop basis $(\ell_*,\ell'_*)$ on $F$ with   $([\ell_i],[\ell'_i])=(z_i,z'_i)\, (i=1,2,\dots,n)$
is constructed by a diffeomorphism realization between symplectic bases of $F$. 
The $Z_2$-quadratic function $q:H_1(F;Z)\to Z_2$ gives 
$q([\ell_i])=\sum_{j=1}^n c_{ij}c'_{ij}=0$, $q([\ell'_i])=\sum_{j=1}^n d_{ij}d'_{ij}=0$ 
for all $i$, showing that the simple loop basis $(\ell_*,\ell')$ is a spin loop basis 
on $F$. 
Let $(E_*\times I ,E'_*\times I )$ be an O2-handle basis of $F$ with
$(\partial E_*,\partial E'_*)=(\ell_*.\ell'_*)$.  
By Lemma~2.1, 
$[S(E_i)]=\sum_{j=1}^n (c_{ij}+2m_{ij})[S(D_j)]+\sum_{j=1}^n (c'_{ij}+2m'_{ij})[S(D'_j)]$, 
$[S(E'_i)]=\sum_{j=1}^n (d_{ij}+2n_{ij})[S(D_j)]+\sum_{j=1}^n (d'_{ij}+2n'_{ij})[S(D'_j)]$
with some integers  $m_{ij}, m'_{ij} , n_{ij}, n'_{ij}$ for all $i,j$, Thus, 
$x_i=[S(E_i)]-2\sum_{j=1}^n m_{ij}[S(D_j)]-2\sum_{j=1}^n m'_{ij}[S(D'_j)]$,  
$x'_i=[S(E'_i)]-2\sum_{j=1}^n n_{ij}[S(D_j)]-2\sum_{j=1}^n n'_{ij}[S(D'_j)]$ 
for all $i$. Since the homology classes $[S(D_j)]$, $[S(D'_j)]$ are $Z$-linear 
combinations of the basis  $([S(E_*)], [S(E'_*)])$ of $H_2(\Sigma;Z)$, 
the desired identities with some $A_i, A'_i, B_i, B'_i$ are obtained. 
This completes the proof of Lemma~2.2.

\phantom{x}

\begin{figure}[hbtp]
	\begin{center}\includegraphics[width=6.5 cm,height=5.5 cm]{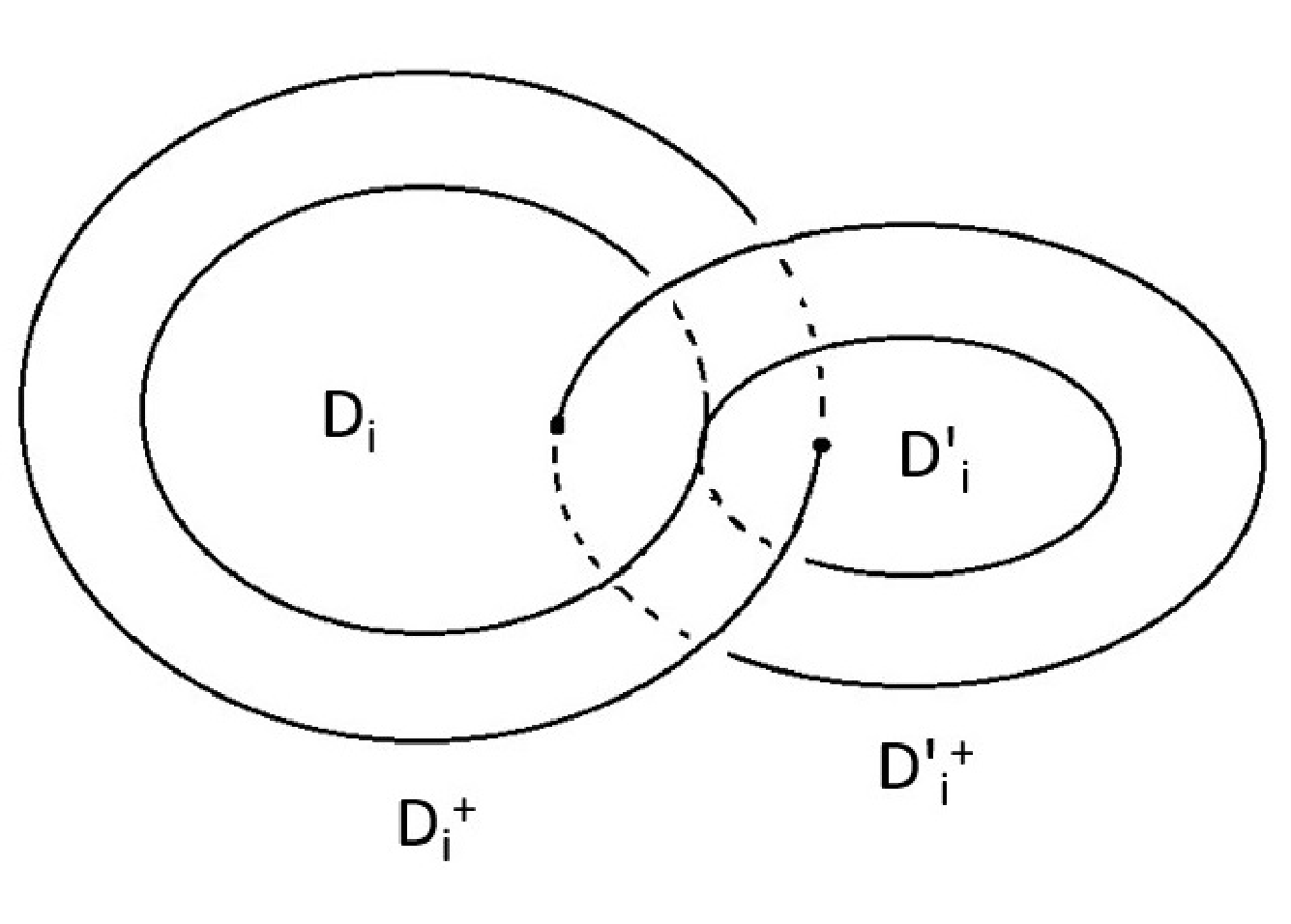}\end{center}
	\caption{Extended disks  $D_i^+$ and ${D'}_i^+$}
	\label{fig2}
\end{figure}

\begin{figure}[hbtp]
	\begin{center}\includegraphics[width=7.5 cm,height=8.5 cm]{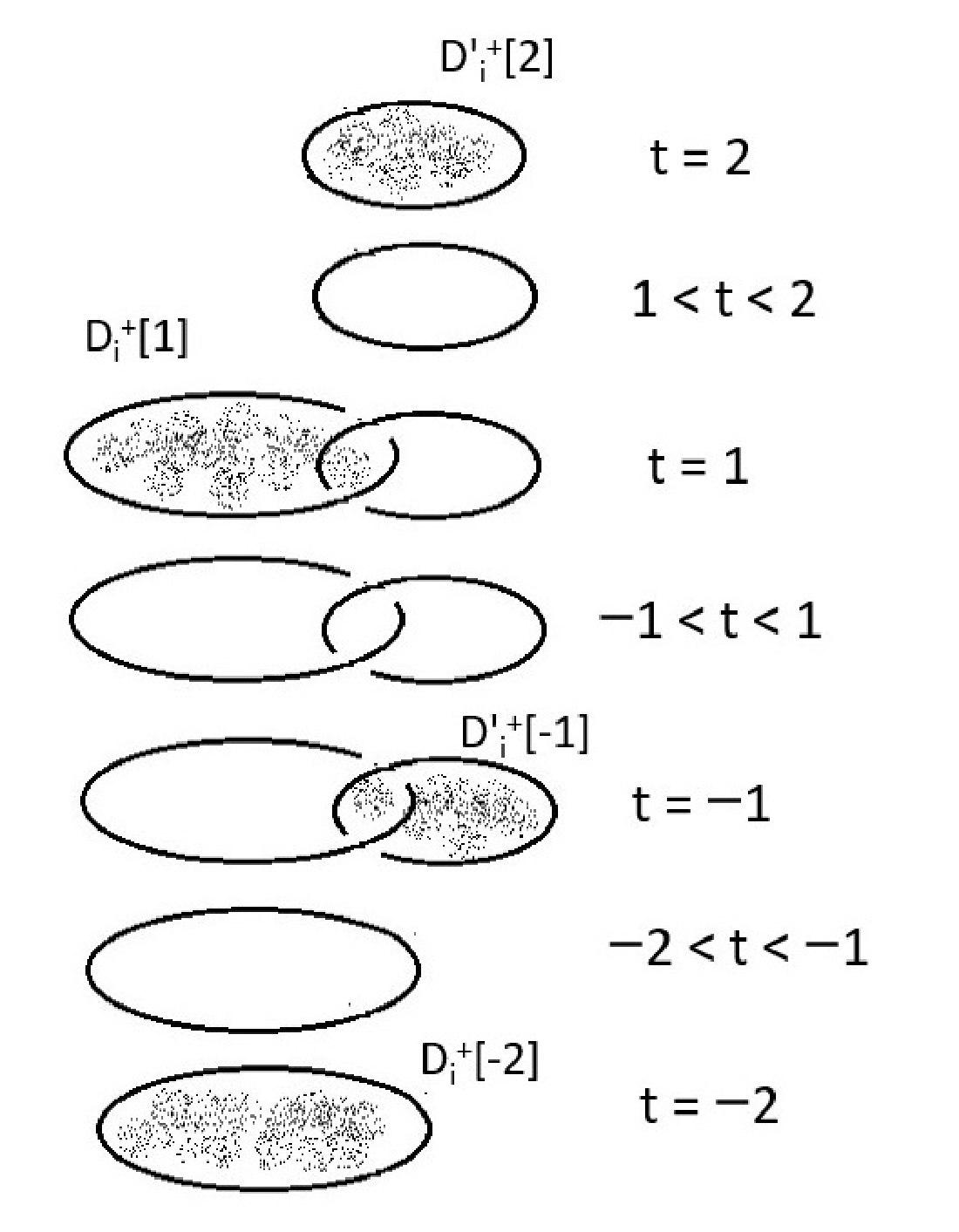}\end{center}
	\caption{Montesinos's twin  $(S_i, S'_i)$}
	\label{fig3}
\end{figure}

Let $S^4$ be the one-point compactification of the 4-space $R^4$. For 
the 3-space $R^3$ and an interval $J\subset R$, the notation  
$R^3J=\{(x,t)\in R^4|\, x\in R^3, t\in J\}$ is used. 
Consider the trivial surface-knot $F$ as  a standard surface in $R^3$ and the O2-handle basis 
$(D_*\times I,D'_*\times I)$ on $F$ is embedded in $R^3$. 
Let $D^+_i$ be a slightly 
 extended disk of $D_i$ so that the boundary loop $\partial D^+_i$  is disjoint 
 from $F$  and 
meets transversely the disk $D'_i$ at a single point, and  ${D'}_i^+$ a slightly extended 
disk of $D'_i$ so that the boundary loop $\partial{D'}_i^+$ is disjoint from $F$  and meets 
 transversely the disk $D_i$ at a single point, Fig.~\ref{fig2}. 
 A {\it 2-sphere surrounding} $D_i$ is a 2-sphere $S_i$ in $R^4$ which is the boundary of the 3-ball $D_i^+[-2,1]$ in $R^4$, and a {\it 2-sphere surrounding} $D'_i$ is a 2-sphere $S'_i$ in $R^4$ which is the boundary of the 3-ball  ${D'}_i^+[-1,2]$ in $R^4$.
The 2-spheres $S_i, S'_i\, (i=1,2,\dots,n)$ are disjoint for distinct indexes $i$ and 
meet transversely with just two points of opposite signs in $R^3[\pm 1]$
for the same index $i$.  The 2-sphere pair $(S_i,S'_i)$ in $S^4$  is called  
{\it Montesinos's twin}, Fig.\ref{fig3}, \cite{7}.  
The connected sum $D_i\# S_i$ is made along an arc $\beta'$ in $S'_i$
joining  the intersection point $D_i\cap S'_i$ in  $R^3[0]$ with  the intersection point of $S_i\cap S'_i$ in $R^3[1]$. 
The 2-sphere $S_i$  is oriented so that  
the disk $D_i\# S_i$ is oriented with the orientation  inherited from $D_i$.  
Similarly, the connected sum $D'_i\# S'_i$ is made along an arc $\beta$ in $S_i$
joining  the intersection point $D'_i\cap S_i$ in $R^3[0]$ with  the intersection point of $S_i\cap S'_i$ in  $R^3[-1]$. 
The 2-sphere $S'_i$  is oriented so that  
the disk $D'_i\# S'_i$ is oriented  with the orientation inherited from $D'_i$.  
To show Lemma~1.2, the following lemma is used.

\phantom{x}

\noindent{\bf Lemma~2.3.} The homology classes 
$[S_i]$ and $[S'_i]$ in $\Sigma$  are given by $[S_i]=-[S(D_i)]$, $[S'_i]=-[S(D'_i)]$ 
for all $i$.

\phantom{x}

\noindent{\bf Proof of Lemma~2.3.} From construction, $S_i$ and $S'_i$ are disjoint 
from $S(D_j), S(D'_j)$ for every $j\ne i$, and 
the intersection numbers $\mbox{Int}([S_i],[S(D_i)]) =\mbox{Int}([S_i],[S_i]) =\mbox{Int}([S'_i],[S'_i]) =\mbox{Int}([S'_i],[S(D'_i)])=0$ and  
$\mbox{Int}([S_i,[S(D'_i)]]) =\varepsilon$, 
$\mbox{Int}([S'_i],[S(D_i)]]) =\varepsilon'$ for $\varepsilon, \varepsilon' =\pm 1$, so that 
$[S_i]=\varepsilon[S(D_i)]$, $[S'_i]=\varepsilon'[S(D'_i)]$. 
Since  $[\alpha S(D_i)]=-[S(D_i)]$ and 
$\mbox{Int}([S(D_i)], [S'_i])=\mbox{Int}([\alpha S(D_i)], [\alpha S'_i])$, 
the identity $[\alpha S'_i]=-[S'_i]$ is obtained. Similarly, 
since  $[\alpha S(D'_i)]=-[S(D'_i)]$ and 
$\mbox{Int}([S(D'_i)], [S_i] )=\mbox{Int}([\alpha S(D'_i)], [\alpha S_i])$, 
the identity $[\alpha S_i]=-[S_i]$ is obtained. 
Let $(E_*\times I ,E'_*\times I )$ be the O2-handle basis on $F$ 
given by $(E_i,E'_i)=(D_i\# S_i,D'_i\# S'_i)$ and $(E_j,E'_j)=(D_j,D'_j)$ for $j\ne i$. 
Then the identities  
$[S(E_i)]=[S(D_i]+ [S_i]-[\alpha S_i]=[S(D_i)]+2[S_i]$,
$[S(E'_i)]=[S(D'_i)]+ [S'_i]-[\alpha S'_i]=[S(D'_i)]+2[S'_i]$
are obtained in $\Sigma$.
Thus, the identities 
 $[S_i]=-[S(D_i)]$, $[S'_i]=-[S(D'_i)]$ in $\Sigma$ 
 are obtained by the identities $[S_i]=\varepsilon[S(D_i)]$, 
 $[S'_i]=\varepsilon'[S(D'_i)]$ and 
$\mbox{Int}([S(E_i)],[S(E'_i)])=\mbox{Int}([S(D_i],[S(D'_i)])=1$. 
This completes the proof of Lemma~2.3.

\phantom{x}

As a by-product of  the proof of Lemma~2.3,  an elementary transformation  
on an O2-handle basis on a trivial surface-knot in $S^4$ was found. 
In fact,  the following elementary transformation  on the O2-handle basis 
$(D_*\times I,D'_*\times I)$ on $F$ is presented.

\phantom{x}

\noindent{\bf Operation I.} An O2-handle basis 
$(D_*\times I,D'_*\times I)$ on $F$ is changed into an O2-handle basis 
$(E_*\times I,E'_*\times I)$ on $F$ with the same attachment as
$(D_*\times I,D'_*\times I)$ such that 
$[S(E_i)]=-[S(D_i)]$, $[S(E'_i)]=-[S(D'_i)]$ in $\Sigma$ and 
$(E_j,E'_j)=(D_j,D'_j)$, $j\ne i,$ for any given $i$.

\phantom{x}

Here are further elementary transformations on the O2-handle basis $(D_*\times I ,D'_*\times I )$ on $F$. 

\phantom{x}

\noindent{\bf Operation II.} Let $i\ne j$ and join the intersection point 
$v_i=D_i\cap S'_i$ 
with the intersection point $v'_j=D'_j\cap S_j$ by an arc $\gamma$  in $S^4$ whose interior is disjoint from $F$ and $(D_*\times I ,D'_*\times I )$, where the 2-spheres  $S'_i$ and $S_j$ are taken to be unoriented. 
Construct the connected sum pairs $(E_i,E'_i)=(D_i\# S_j,  D'_i)$ and 
$(E_j,E'_j)=(D_j, D'_j\# S'_i)$ along the arc $\gamma$. 
For $h\ne i, j$, let $(E_h,E'_h)=(D_h,D'_h)$. Then 
$(E_*\times I ,E'_*\times I )$ is an O2-handle basis on $F$ with the same attachment as $(D_*\times I ,D'_*\times I )$ such that 
$[S(E_i)]= [S(D_i]-2\varepsilon[S(D_j)]$, $[S(E'_j)]= [S(D'_j)]+ 2\varepsilon[S(D'_i)]$ 
with $\varepsilon=\pm 1$ in 
$\Sigma$  and    
$E'_i=D'_i$, $E_j=D_j$,  $(E_h,E'_h)=(D_h,D'_h)$, $h\ne i, j$, 
for any given distinct indexes $i$ and $j$.

In fact, since by Lemma~2.3 
$[S(E_i)]=[S(D_i]+2\varepsilon [S_j] =[S(D_i)]-2\varepsilon [S(D_j)]$ and 
$ [S(E'_j)])=[S(D'_j)]+2\varepsilon' [S'_i]=[S(D'_j)]-2\varepsilon'[S(D'_i)]$ 
for signs  $\varepsilon, \varepsilon'$ being $\pm1$,  the intersection number 
$\mbox{Int}([S(E_i)], [S(E'_j)])=0$ implies $\varepsilon+\varepsilon' =0$. 

\phantom{x}

\noindent{\bf Operation III.} Let $i\ne j$ and join the intersection point 
$v_i=D_i\cap S'_i$ with the intersection point $v'_j=D_j\cap S'_j$ by an arc $\gamma$  in $S^4$ whose interior is disjoint from $F$ and $(D_*\times I ,D'_*\times I )$, 
where the 2-spheres  $S'_i$ and $S'_j$ are taken to be unoriented. 
Construct the connected sum pairs $(E_i,E'_i)=(D_i\# S'_j,  D'_i)$ and 
$(E_j,E'_j)=(D_j\# S'_i, D'_j)$ along the arc $\gamma$. 
For $h\ne i, j$, let $(E_h,E'_h)=(D_h,D'_h)$. Then 
$(E_*\times I ,E'_*\times I )$ is an O2-handle basis on $F$ with the same attachment as $(D_*\times I ,D'_*\times I )$ such that 
$[S(E_i)]=[S(D_i)]- 2\varepsilon[S(D'_j)]$, $[S(E_j)]=[S(D_j)]+ 2\varepsilon[S(D'_i)]$ 
with $\varepsilon=\pm 1$ in  $\Sigma$     and 
$E'_i=D'_i$, $E'_j=D'_j$,  $(E_h,E'_h)=(D_h,D'_h)$, $h\ne i, j$,
for any given distinct indexes $i$ and $j$.

In fact, since by Lemma~2.3 
$[S(E_i)]=[S(D_i)]+2\varepsilon [S'_j]=[S(D_i)]-2\varepsilon [S(D'_j)]$ and 
$ [S(E_j)])=[S(D_j)]+2\varepsilon' [S'_i]=[S(D_j)]-2\varepsilon'[S(D'_i)]$ 
for signs $\varepsilon, \varepsilon' $ being $\pm1$, the intersection number 
$\mbox{Int}([S(E_i)], [S(E_j)])=0$ implies $\varepsilon+\varepsilon' =0$. 
 
\phantom{x}
 
\noindent{\bf Operation IV.} Let $i\ne j$ and join the intersection point 
$v_i=D'_i\cap S_i$ 
with the intersection point $v'_j=D'_j\cap S_j$ by an arc $\gamma$  in $S^4$ whose interior is disjoint from $F$ and $(D_*\times I ,D'_*\times I )$, where the 2-spheres  $S_i$ and $S_j$ are taken to be unoriented. 
Construct the connected sum pairs $(E_i,E'_i)=(D_i,  D'_i\# S_j)$ and 
$(E_j,E'_j)=(D_j, D'_j\# S_i)$ along the arc $\gamma$. 
For $h\ne i, j$, let $(E_h,E'_h)=(D_h,D'_h)$. Then 
$(E_*\times I ,E'_*\times I )$ is an O2-handle basis on $F$ with the same attachment as $(D_*\times I ,D'_*\times I )$ such that 
$[S(E'_i)]=[S(D'_i)]- 2\varepsilon[S(D_j)]$, $[S(E'_j)]=[S(D'_j)]+2\varepsilon[S(D_i)]$ 
with $\varepsilon=\pm1$ in $\Sigma$   and   
$E_i=D_i$, $E_j=D_j$,  $(E_h,E'_h)=(D_h,D'_h)$, $h\ne i, j$, 
for any given distinct indexes $i$ and $j$.

In fact, since by Lemma~2.3 
$[S(E'_i)]=[S(D'_i)]+2\varepsilon [S_j] = [S(D'_i)]-2\varepsilon [S(D_j)]$ and 
$ [S(E'_j)])=[S(D'_j)]+2\varepsilon' [S_i]=[S(D'_j)]-2\varepsilon'[S(D_i)]$ 
for signs $\varepsilon, \varepsilon' $ being $\pm1$, the intersection number 
$\mbox{Int}([S(E'_i)], [S(E'_j)])=0$ implies $\varepsilon+\varepsilon' =0$. 
 
\phantom{x}

In Operations II, III, IV, the sign $\varepsilon$ can take both $+1$ and $-1$, because
a normal 3-disk bundle of $\gamma$ in $S^4$ admits the Hopf link bundle over  
$\gamma$ used for the connected sums and the Hopf link has a component-preserving inversion. 
The following observation on an odd integer $u$ and a non-zero integer $c$ is used for the proof of Lemma 2.1.

\phantom{x}

\noindent{\bf Lemma 2.4.}  Let $u$ be  an odd integer,  $c$ a non-zero integer, 
and $\varepsilon  = \pm1$ the sign of the product integer $uc$.  
If $|c|\geq  |u|$, then $0\leq |2c-2\varepsilon u|<2|c|$.  
If $|u|> |c|$, then the integer $u = 2\varepsilon c$ is odd with 
$1 \leq  |u = 2\varepsilon c|< |u|$.

\phantom{x}

In fact, let  $\varepsilon'c> 0$ for $\varepsilon'=\pm 1$. If $\varepsilon'c \geq  \varepsilon'\varepsilon u>0$, then  
$0 \leq 2\varepsilon'c = 2\varepsilon'\varepsilon u < 2\varepsilon'c$, showing that 
$|c|\geq  |u|$ implies $0 \leq |2c = 2\varepsilon u|< 2|c|$.  
If  $\varepsilon'\varepsilon u > 2\varepsilon'c >0$, then 
$0<\varepsilon'\varepsilon u =2\varepsilon'c < \varepsilon'\varepsilon u$. If $2\varepsilon'c > \varepsilon'\varepsilon u >\varepsilon'c >0$, then 
$0 >\varepsilon'\varepsilon u =2\varepsilon'c > =\varepsilon'c> 
=\varepsilon'\varepsilon u$. These inequalities imply 
$1 \leq |u = 2\varepsilon c|< |u|$, as desired.  

In applications of Observation~2.4 to the proof of Lemma~1.2, it is used that the coefficients $u_i = 1 + 2a_{ii}$ of $[S(E_i)]$ in $x_i$ and 
$u'_i = 1 + 2a'_{ii}$ of $[S(E'i)]$ in $x'_i$ are odd integers. The proof of Lemma~1.2 is done as follows.

\phantom{x}

\noindent{\bf Proof of Lemma~1.2.}  Fix the orthogonal basis $(x_*,x'_*)$  of $\Sigma$.  
By Lemma~2.2, there is an O2-handle basis  $(E_*\times I , E'_*\times I )$ on  $F$  
such that 
$x_i=[S(E_i)]+ 2A_i+2A'_i$ and $x'_i=[S(E'_i)]+2B_i+2B'_i$ for all $i$, where 
$A_i=\sum_{j=1}^n a_{ij}[S(E_j)]$, $A'_i=\sum_{j=1}^n a'_{ij}[S(E'_j)]$,  
$B_i=\sum_{j=1}^n b_{ij}[S(E_j)]$,  
$B'_i=\sum_{j=1}^n b'_{ij}[S(E'_j)]$ 
with some integers $a_{ij}$, $a'_{ij}$, $b_{ij}$, $b'_{ij}$ for all $i, j$. 

Let $n=1$. Since $(x_1,x'_1)$ and $([S(E_1)],[S(E'_1)])$ are orthogonal basis of 
$\Sigma$, it holds that $a'_{11}=b'_{11}=0$ and either 
$a_{11}=b_{11}=0$ or $a_{11} = b_{11} = -1$, so that  either 
$x_i = [S(E_i)]$,    $x'_i = [S(E'_i)]$ or $x_i = - [S(E_i)]$, $x'_i = - [S(E'_i)]$. 
For the latter case, apply Operation I to obtain $x_1 = [S(E_1)]$ and 
$x'_1 = [S(E'_1)]$.
Let $n\geq 2$.
By  Lemma~2.4,  By a finite number of applications of Operation II and Observation~2.4 using the transformations $[S(E_n)]\to [S(E_n)]\pm 2[S(E_j)]$ and 
$[S(E_j)]\to [S(E_i)]\pm 2[S(E_n)]$ for $j<n$, 
there is an O2-handle basis  $(E_*\times I , E'_*\times I )$ on  $F$  such that 
$x_n=\varepsilon [S(E_n)]+2A'_n\, (\varepsilon=\pm1)$, $x'_n=[S(E'_n)]+2B_n+2B'_n$ with some $A'_n, B_n, B'_n$.  
By a finite number of applications of Operation III and Observation~2.4 using 
the transformations $[S(E_n)]\to [S(E_n)]\pm 2[S(E'_j)]$ for $j<n$, 
there is an identity $x_n = \varepsilon [S(En)] + 2a'_{nn}[S(E'_n)]$ 
for some integer $a'_{nn}$. 
Then $a'_{nn} = 0$, because $\mbox{Int}(x_n, x_n)=\mbox{Int}([S(E_n)], [S(E_n)])=
\mbox{Int}([S(E'_n)], [S(E'_n)])=0$ and  $\mbox{Int}([S(E_n)], [S(E'_n)])=1$.
Thus, $x_n = \varepsilon  [S(E_n)]$, 
$x'_n = [S(E'_n)] + 2B_n + 2B'_n$ with some $B_n$, $B'_n$. 
By a finite number of Operation IV and Observation~2.4 using 
$[S(E'_n)]\to [S(E'_n)]\pm 2[S(E_j)]$ for $j<n$,
the identities  
$x_n = \varepsilon  [S(En)]$, $x'_n =  [S(E'_n)] + 2b_{nn}[S(En)] + 2B'_n$ 
are obtained for some integer $b_{nn}$ and some $B'_n$.
By Operation I, assume that $\varepsilon=1$. Use that 
$(x_*, x'_*)$ and $([S(E_*)], [S(E'_*)])$ are orthogonal bases of $\Sigma$. Since 
$\mbox{Int}_(x'_n, x'_n) = 4b_{nn} = 0$ and thus $b_{nn} = 0$, the identities 
$x_n = [S(E_n)]$,  $x'_n = [S(E'_n)] + 2B'_n$ are obtained. Then 
$B'_n = 0$, for $\mbox{Int}(x_j, x'_n) = 0\, (j<n)$ or $1\,(j=n)$, and 
the identities $x_n = [S(E_n)]$,  $x'_n =  [S(E'_n)]$ are obtained. 

Let $F_{n-1}$ be the trivial surface-knot of genus $n-1$ in $S^4$ obtained from 
$F$ by surgery along the O2-handle pair $(E_n\times I ,E'_n\times I )$, \cite{3}. 
Then the system $(x_i,x'_i)\,(i=1,2,\dots,n-1)$ are regarded as an orthogonal 
basis for the stable 4-sphere $\Sigma(n-1)=S^4(F_{n-1})_2$ of genus $n-1$. 
By inductive assumption on $n$, there is a replacement of the O2-handle basis 
$(E_i\times I ,E'_i\times I )\, (i=1,2,\dots,n-1)$ on $F_{n-1}$ keeping the attaching part fixed so that 
$x_i=[S(E_i)]$ and $x'_i= [S(E'_i)]$  for all $i\, (i=1,2,\dots,n-1)$ in $\Sigma(n-1)$.  
Thus, there is  an O2-handle basis $(E_*\times I ,E'_*\times I )$ on $F$ 
such that $x_i=[S(E_i)]$ and $x'_i= [S(E'_i)]$ for all $i\, (i=1,2,\dots,n)$ in 
$\Sigma=\Sigma(n)$, \cite{3}. This completes the proof of Lemma~1.2. 

\phantom{x} 

This completes the proof of Theorem~1.1.

\phantom{x}

\noindent{\bf 3. Application to diffeomorphism of stable 4-sphere}
 
The orientation-preserving diffeomorphism group of the 4-ball $D^4$ keeping the 
boundary $\partial D^4$ point-wise fixed is denoted by 
$\mbox{Diff}^+(D^4, \mbox{rel}\partial)$. An identity-shift of the stable 4-sphere $\Sigma$     is a diffeomorphism  $\iota: \Sigma   \to  \Sigma$     obtained from the identity 
$1 :  \Sigma    \to \Sigma$     by replacing the identity on a 4-ball in 
$ \Sigma=S^4(F)_2$   disjoint from $F$  with an element of 
$\mbox{Diff}^+(D^4, \mbox{rel}\partial)$. The following theorem is shown as an application of Theorem 1.1.

\phantom{x}

\noindent{\bf Theorem 3.1.} For every stable 4-sphere $\Sigma$, every orientation-preserving diffeomorphism $f$   of $ \Sigma$     is characterized by the induced orthogonal base change automorphism of $H_2(\Sigma; Z)$, which  is nothing but the lift  $g'$ of an equivalence $g$ of the trivial surface-knot space $(S^4, F)$ to 
$\Sigma$     up to  composite of an identity-shift  $\iota$   and smooth isotopy on 
$S^4(F)_2 = \Sigma$.

\phantom{x}

\noindent{\bf Proof of Theorem 3.1.}  For an orthogonal basis  $(x_* , x'_*)$ of the stable 4-sphere $ \Sigma$     of genus  n, let $(y_* , y'_* )$ be the orthogonal basis of $ \Sigma$     given by $f_*(x_i  ) = y_i$ and $f_*(x'_i) = y'_i$ for all $i$. By  Theorem 1.1, there is an equivalence $g$ of a trivial surface-knot space $(S^4, F)$ of genus  $n$ in $S^4$  whose lift  $g'$   to   $S^4(F)_2 = \Sigma$     has  $g'_*(x_i  ) = y_i$ and  
$g'_*(x'_i) = y'_i$ for all $i$. The diffeomorphisms $f$  and  $g'$   are homotopic since $ \Sigma$     is a simply connected closed oriented 4-manifold, and $f$  and  $g'$   induce the same intersection form on the second homology group
$H_2(\Sigma; Z)$, \cite{8}. Then there is an identity shift  $\iota$  of 
$ \Sigma=S^4(F)_2$   such that the composite $f \iota$  is smoothly isotopic to $g'$, [1, Theorem 4.1]. This completes the proof of Theorem 3.1.

In Piecewise-linear category, every piecewise-linear self-homeomorphism of the 
4-disk keeping the boundary identically is piecewise-linearly $\partial$-relatively isotopic to the identity, which is well-known as Alexander trick, as noted \cite{1}. Thus, the following corollary is obtained. 

\phantom{x}

\noindent{\bf Corollary 3.2.} For every stable 4-sphere $\Sigma$, every 
orientation-preserving piecewise-linear self-homeomorphism $f$  of $ \Sigma$     
is characterized by the induced orthogonal base change automorphism  of 
$H_2(\Sigma; Z)$, which is nothing but the lift  $g'$   of an equivalence $g$ of the trivial surface-knot space  $(S^4, F)$ to $\Sigma$  up piecewise-linear isotopy of 
$\Sigma$.

\phantom{x} 

\noindent{\bf 4. Application to TOP stable 4-sphere}
 
A TOP stable 4-sphere of genus $n$ is a topological 4-manifold $X = X(n)$ which is homeomorphic to the stable 4-sphere $\Sigma(n)$. An orthogonal basis  $(x_* , x'_*)$ of $X(n)$ is similarly  defined by the intersection form on $X(n)$ and sent to an orthogonal basis of $\Sigma(n)$ by the  homeomorphism. A TOP trivial surface-knot space of genus $n$ is a topological pair  $(Y, F^Y)$ homeomorphic to a trivial surface-knot space $(S^4, F)$ of genus n. Then $Y$ is a TOP 4-sphere and $F^Y$ is a TOP surface-knot in $Y$ which bounds a TOP handlebody in $Y$, the pullback of a handlebody bounded by $F$ in $S^4$. A  TOP O2-handle pair on $F^Y$ is the pullback of an O2-handle pair on the trivial surface-knot $F$ in $S^4$. Fix the orientation of the 4-sphere  $(Y, F^Y)$  inherited from the orientation of $(S^4,  F)$. The following theorem is a TOP version of Theorem 1.1. 

\phantom{x} 

\noindent{\bf Theorem 4.1.} For every TOP stable 4-sphere $X(n)$, there is a TOP trivial surface-knot space  $(Y, F^Y)$  of genus $n$ whose double branched covering space  $Y(F^Y)_2$ = $X(n)$. Even if X is a  smooth 4-manifold, unless $X(n)$ is not diffeomorphic to the stable 4-sphere $\Sigma(n)$, the TOP trivial surface-knot space  $(Y, F^Y)$  is not smooth-able. 
Every orthogonal basis  $(x_* , x'_*)$ of $X(n)$ is  represented by a TOP O2-sphere basis $(S(E_*^{Yx}),  S({E'}_*^{Yx}))$  of $X(n)$ constructed from a TOP O2-handle basis $(E_* Yx \times I, {E'}_*^{Yx} \times I )$ on $F^Y$ in $Y$. 
For any two orthogonal bases  $(x_* , x'_*)$ and $(y_* , y'_* )$  of $X(n)$, there are TOP O2 -handle bases $(E_*^{Yx} \times I, {E'}_*^{Yx} \times I)$ and 
$(E_*^{Yy} \times I, {E'}_*^{Yy} \times I )$ on $F^Y$ in  $Y$ which are transformed into each other by a TOP equivalence $g$ of  $(Y, F^Y)$  such that  $(x_* , x'_*)$  and $(y_* , y'_* )$ are represented by the TOP O2-sphere bases 
$(S(E_*^{Yx}),  S({E'}_*^{Yx}))$  and $(S(E_*^{Yy}),  S({E'}_*^{Yy}))$, respectively which are transformed into each other by the lift $g'$ of $g$ to  $Y(F^Y)_2$ = $X(n)$.

\phantom{x} 

\noindent{\bf Proof of Theorem 4.1.} Since $X(n)$ is a TOP stable 4-sphere of genus n, there is an  orientation-preserving homeomorphism $h: X(n)\to  \Sigma(n)$. Let $(S^4, F)$ be a trivial surface-knot space of genus $n$ with   $S^4(F)_2 =\Sigma(n)$. For the nontrivial covering involution $\alpha$  of   $S^4(F)_2$, let 
$\alpha^h = h^{-1}\alpha h$ be a TOP involution on $X(n)$ with the fixed-point set $h^{-1}(F)$.  Let $Y$ be the orbit space of $X$ by $\alpha^h$, admitting a homeomorphism $h^Y: Y \to  S^4$  induced  from $h$ by regarding $S^4$  as the orbit space of   $S^4(F)_2$   by $\alpha$. The projection $p^Y: X(n) \to  Y$  defines the double branched covering projection of $Y$ branched along the projection  image surface $F^Y$ of $h^{-1}(F)$, and $X(n)$ = $Y(F^Y)_2$. In fact, the homeomorphism 
$h: X(n) \to  \Sigma(n)$  defines a homeomorphism 
$h^Y:  (Y, F^Y)  \to  (S^4, F)$ and is considered as the lift  $Y(F^Y)_2 \to S^4(F)_2$   of $h^Y$. The pair  $(Y, F^Y)$  is a TOP trivial surface-knot space of genus $n$ sent by $h^Y$ from the trivial surface-knot space $(S^4, F)$ of genus $n$. For any orthogonal basis  $(x_*, x'_*)$ of  $X(n)$, let $(x_*^S, {x'}_*^S)$ be the orthogonal basis of $\Sigma(n)$ given by $x_i^S = h_*(x_i)$, ${x'}_i^S = h_*(x'_i)\,
(i=1,2, \dots, n)$. By Theorem 1.1, $(x_*^S, {x'}_*^S)$ is represented by the O2-sphere basis $(S(E_*^x),S({E'}_*^x))$ of $\Sigma(n)$ constructed from an O2-handle basis 
$(E_*^x\times  I, {E'}_*^x \times  I )$  on $F$  in $S^4$. Then $(x_*, x'_*)$ is represented by the TOP O2-sphere basis $(h^{-1}S(E_*^x), h^{-1}S({E'}_*^x))$ of $X(n)$ constructed from the TOP O2 handle basis 
$((h^Y)^{-1}E_*^x\times  I, (h^Y)^{-1}{E'}_*^x \times  I )$ on $F^Y$ in $Y$. For another orthogonal basis $(y_* , y'_* )$ of $X(n)$, let $(y_*^S, {y'}_*^S)$ be the orthogonal basis of $\Sigma(n)$  given by $y_i^S = h_*(y_i)$, 
${y'}_i^S = h_*(y'_i)\, (i=1,2, \dots, n)$, which is represented by the O2-sphere  basis 
$(S(E_*^y), S({E'}_*^y))$ of $\Sigma(n)$ constructed from an O2-handle basis 
$(E_*^y \times  I, {E'}_*^y \times  I )$ on  $F$  in $S^4$. Then there is an an equivalence $f$  of $(S^4, F)$ sending  $(E_*^x \times  I,  {E'}_*^x \times  I )$  to 
 $(E_*^y \times I, {E'}_*^y \times  I)$  whose lift $f'$ to   $S^4(F)_2 =\Sigma(n)$ sends $(S(E_*^x), S({E'}_*^x))$ to $(S(E_*^y), S({E'}_*^y))$. Then  the composite TOP equivalence $g = (h^Y)^{-1}f h^Y$ of  $(Y, F^Y)$  sends the TOP O2 -handle basis  
$((h^Y)^{-1}E_*^x\times  I, (h^Y)^{-1}{E'}_*^x \times  I )$ on $F^Y$ in $Y$ to the TOP O2-handle basis 
$((h^Y)^{-1}E_*^y\times  I, (h^Y)^{-1}{E'}_*^y  \times  I)$ on $F^Y$ in $Y$. The lift 
$g' = h^{-1}fh$ of $g$ to $X(n)$ sends the TOP O2-sphere basis 
$(h^{-1}S(E_*^x), h^{-1}S({E'}_*^x))$ of $X(n)$ to the TOP O2-sphere basis 
$(h^{-1}S(E_*^y), h^{-1}S({E'}_*^y))$ of $X(n)$ and thus  sends $(x_*, x'_*)$ to 
$(y_*, y'_*)$. Suppose that $X(n)$ is a smooth 
4-manifold obtained from a  smooth surface-knot space  $(Y, F^Y)$, but not diffeomorphic to $\Sigma(n)$. Then there is an orientation-preserving diffeomorphism $g$ from $Y$ to the 4-sphere  $S^4$, [1]. Since the image $G=gF^Y$ is a smooth surface-knot in $S^4$  with fundamental group $\pi_1(S^4\setminus G, x_0)$ an infinite cyclic group, the surface-knot $G$ is a trivial surface-knot of genus $n$, and $G$ is  equivalent to $F$, \cite{3}, \cite{6}. Thus, $X(n)$ is diffeomorphic to $\Sigma(n)$, a contradiction. Thus, the  TOP surface-knot space  $(Y, F^Y)$  cannot be smooth. By writing $((h^Y)^{-1}E_*^x\times  I, (h^Y)^{-1}{E'}_*^x \times I )$, 
$((h^Y)^{-1}E_*^y\times  I, (h^Y)^{-1}{E'}_*^y \times  I )$, 
$(h^{-1}S(E_*^x), h^{-1}S({E'}_*^x)), (h^{-1}S(E_*^y), h^{-1}S({E'}_*^y))$ as 
$(E_*^{Yx}\times  I, {E'}_*^{Yx} \times  I )$, $(E_*^{Yy}\times  I, {E'}_*^{Yy} \times  I)$, 
$(S(E_*^{Yx}), S({E'}_*^{Yx}))$, $(S(E_*^{Yy}), S({E'}_*^{Yy}))$, respectively, 
the proof of  Theorem 4.1 is complete.

\phantom{x} 

It is known that there are smooth stable 4-spheres $X(n)$, not diffeomorphic to 
$\Sigma(n)$ for  infinitely-many $n$, \cite{9}. Also, it is known that for every TOP stable 4-sphere $X(n)$,  any two homotopic self-homeomorphisms of $X(n)$ are isotopic by TOP isotopy of $X(n)$,  \cite{10}. Thus, every orientation-preserving 
self-homeomorphism of $X(n)$ is characterized  by the induced orthogonal base change automorphism of $H_2(X(n); Z)$, which is nothing but the lift 
$g': X(n) \to  X(n)$ of a TOP equivalence $g$ of the TOP trivial surface-knot  
$(Y, F^Y)$  of genus $n$ up to TOP isotopy.

\phantom{x} 

\noindent{\bf Conclusions.}
For every orthogonal basis  $(x_* , x'_*)$ of the stable 4-sphere $\Sigma(n)$ is realized by the orthogonal basis $([S(E_*)], [S(E'_*)])$  for an O2-handle basis 
$(E_*\times I, E'_*\times I)$ on a trivial surfaceknot space $(S^4, F)$ of genus $n$ whose double branched covering space   $S^4(F)_2$   is $\Sigma(n)$ by Lemma 1.2. Then it is used that any two O2-handle bases on $F$  can be transformed into each other by an equivalence of the surface-knot space  $(S^4, F)$, \cite{6}. It is concluded  that any two orthogonal bases of $\Sigma$     are transformable into each other by the lift $f'$ of $f$  to $\Sigma =  S^4(F)_2$   (Theorem 1.1). Any two homotopic diffeomorphisms of $ \Sigma$     are smoothly  isotopic if one diffeomorphism is replaced by the composition with an identity-shift,  \cite[Theorem 4.1]{1}. For this proof, Gabai's 4D light bulb theorem is used, \cite{8}. By 
combining this result with Theorem 1.1, it is concluded that every orientation-preserving  diffeomorphism of $ \Sigma$     is nothing but he lift  $g'$   of an equivalence $g$ of the trivial surface-knot space $(S^4, F)$ to $\Sigma =  S^4(F)_2$   up to composition of an identity-shift and smooth isotopy  (Theorem 3.1). Here the question of whether every identity-shift is smoothly isotopic  to the identity remains an unsolved problem, \cite{1}. In Piecewise-linear category and TOP category, this problem is not needed and the arguments proceed well except that  even if a TOP stable 4-sphere $X(n)$ is a smooth 4-manifold, unless $X(n)$ is diffeomorphic to the stable 4-sphere $\Sigma(n)$, the TOP trivial surface-knot space  $(Y, F^Y)$  of genus $n$ with  $Y(F^Y)_2 = X$ cannot be smooth (Corollary 3.2 and Theorem 4.1).

\phantom{x} 

\noindent{\bf Acknowledgements.} 
This work was partly supported by JSPS KAKENHI Grant Number JP21H00978 and MEXT Promotion of Distinctive Joint Research Center Program JPMXP0723833165 and Osaka Metropolitan University Strategic Research Promotion Project (Development of International Research Hubs).

\phantom{x}


\begin{thebibliography}{99}

\bibitem{1} A. Kawauchi, Smooth homotopy 4-sphere, WSEAS Transactions on Mathematics 22 (2023), 690-701.

\bibitem{2} C. T. C. Wall, Diffeomorphisms of 4-manifolds, J. London Math. Soc. 39 (1964),  131-140.

\bibitem{3} A. Kawauchi, Ribbonness of a stable-ribbon surface-link, I. 
A stably trivial surface-link, Topology and its Applications 301(2021), 107522 (16pages). 

\bibitem{4} A. Kawauchi, Smooth homotopy 4-sphere (research announcement), 2191 Intelligence of Low Dimensional Topology, RIMS Kokyuroku 2191 (July 2021), 1-13. 


\bibitem{5}  S. Hirose, On diffeomorphisms over surfaces trivially 
embedded in the 4-sphere, Algebraic and Geometric Topology 2(2002), 791-824. 

\bibitem{6} A. Kawauchi, Uniqueness of an orthogonal 2-handle pair on a surface-link, Contemporary Mathematics (UWP) 4 (2023), 182-188. 

\bibitem{7} J. M. Montesinos, On twins in the four-sphere I, 
Quart. J. Math. Oxford (2), 34 (1983), 171-199.

\bibitem{8} D. Gabai, The 4-dimensional light bulb theorem, J. Amer. Math. Soc. 33 (2020), 609-652.

\bibitem{9} 
A. Akhmedov and B. D. Park, Geography of simply connected spin symplectic 
4-manifolds, Mathematical Research Letters, 17 (2010), 483-492. 

\bibitem{10} 
S. Friedl,  M. Nagel,  P. Orson,  and M. Powell, The foundations of four-manifold theory in the topological category. New York Journal of Mathematics Monographs, 6 (2025), 1-152.


\end{thebibliography}
\end{document}